\documentclass[preprint,12pt]{elsarticle}

\usepackage{amssymb}
\usepackage{xcolor}
\usepackage{natbib}
\usepackage{endnotes}
\usepackage{graphicx}
\usepackage{subfigure}
\usepackage{enumerate}
\usepackage{bm}
\usepackage{amsmath}
\usepackage[linesnumbered,ruled,vlined]{algorithm2e}
\usepackage{amsmath, amsthm, amssymb}

\newtheorem{fact}{Fact}
\usepackage{multirow}
\usepackage{booktabs}
\newcommand{\MC}[1]{\textcolor{black}{ #1}}
\newcommand{\PL}[1]{\textcolor{black}{ #1}}
\def\win{\textit{\#win}}
\def\feas{\textit{\#feas}}
\def\pt{\textit{P(T)}}
\def\improve{\textit{Improve}}
\def\search{\textit{Search}}
\def\restore{\textit{Restore}}
\def\tm{$tm$}
\def\lm{$lm$}
\def\op{$op$}
\usepackage{amsthm}
\newtheorem{myDef}{Definition}

\newcommand{\bl}[1]{\bm{#1}}

\newtheorem{proposition}{Proposition}
\newtheorem{remark}{Remark}

\journal{Artificial Intelligence}

\begin{document}

\begin{frontmatter}



\title{New Characterizations and Efficient Local Search for General Integer Linear Programming}

\author[label1,label2]{Peng Lin}
\author[label1,label2]{Shaowei Cai\corref{cor1}}
\ead{caisw@ios.ac.cn}
\author[label1]{Mengchuan Zou}
\author[label3]{Jinkun Lin}
\cortext[cor1]{Corresponding author}
\affiliation[label1]{organization={State Key Laboratory of Computer Science, Institute of Software, Chinese Academy of Sciences},
            city={Beijing},
            country={China}}

\affiliation[label2]{organization={School of Computer Science and Technology, University of Chinese Academy of Sciences},
            city={Beijing},
            country={China}}
\affiliation[label3]{organization={SeedMath Technology Limited},
            city={Beijing},
            country={China}}



\begin{abstract}
Integer linear programming (ILP) models a wide range of practical combinatorial optimization problems and significantly impacts industry and management sectors. 
This work proposes new characterizations of ILP with the concept of boundary solutions. 
\PL{Motivated by the new characterizations, we develop a new local search algorithm Local-ILP, which is efficient for solving general ILP validated on a large heterogeneous problem dataset.}
We propose a new local search framework that switches between three modes, namely \search{}, \improve{}, and \restore{} modes. {Two new operators are proposed, namely the \textit{tight move} and the \textit{lift move} operators, which are associated with appropriate scoring functions.} Different modes apply different operators to realize different search strategies and the algorithm switches between three modes according to the current search state.
Putting these together, we develop a local search \PL{ILP solver} called Local-ILP.
Experiments conducted on the MIPLIB dataset show the effectiveness of our algorithm in solving large-scale hard ILP problems. {In the aspect of finding a good feasible solution quickly, Local-ILP is competitive and complementary to the state-of-the-art commercial solver Gurobi} and significantly outperforms the state-of-the-art non-commercial solver SCIP. Moreover, our algorithm establishes new records for 6 MIPLIB open instances.
The theoretical analysis of our algorithm is also presented, which shows our algorithm could avoid visiting unnecessary regions.
\end{abstract}




\end{frontmatter}


\section{Introduction}
\MC{
Integer linear programming (ILP) is a fundamental problem in combinatorial optimization, where the objective is to optimize a linear objective function under linear constraints, while variables must take integer values. The ILP is NP-hard, and one of its special cases, the 0-1 ILP problem belongs to Karp's 21 NP-complete problems list. Moreover, ILP has strong descriptive capability and is usually more convenient for describing practical combinatorial optimization problems than the SAT, MaxSAT, 0-1 ILP (PBO) formulations. Many combinatorial optimization problems, such as the knapsack problem, traveling salesman problem, warehouse location problem, decreasing costs and machinery selection problem, network and graph problems (maximum flow, set covering, matching, weighted matching, spanning trees, etc.), and scheduling problems, can all be formulated and solved as ILP problems \citep{wolsey2020integer}.
}


Along with the powerful ability of ILP to model combinatorial optimization problems, there have also been plenty of efforts devoted to solving ILP problems. As the general ILP problem is NP-hard \citep{karp1972reducibility}, there is no known polynomial algorithm that could solve ILP exactly. \MC{As a result, solving methods have been categorized into complete and incomplete methods: complete methods compute the exact optimal solution and prove the optimality; incomplete methods attempt to get a good solution quickly and do not need to guarantee the optimality. }

\PL{Existing ILP solvers are primarily based on complete methods.}
The best-known \MC{complete} approach is branch-and-bound (BnB), which is a classical approach for solving combinatorial optimization problems by iteratively dividing the feasible region and bounding the objective function \citep{land1960automatic,lawler1966branch}. 
In addition, other methods such as cutting plane \citep{gomory1958outline} and auxiliary procedures such as domain propagation \citep{brearley1975analysis, savelsbergh1994preprocessing, achterberg2009scip}, have been proposed for ILP, which are often integrated into the BnB process.
\PL{This is what forms the hybrid approach, which are promising because they leverage the strengths of different methods in a complementary mode. 
Hybrid approaches are the infrastructure of state-of-the-art commercial and non-commercial ILP solvers such as Gurobi \cite{gurobi2022gurobi} and SCIP \citep{achterberg2009scip}, both of which are based on the BnB algorithm, combining, cutting planes, domain propagation, etc.}
However, due to the worst-case exponential running time, complete algorithms are impractical for large-scale instances. 



\PL{Local search is an important class of incomplete algorithms and is powerful for solving NP-hard combinatorial problems in many areas of computer science and operations research \citep{hoos2004stochastic}.}
Local search algorithms iteratively explore the neighborhood of a solution and move towards a better one. It has shown great success in solving SAT and MaxSAT \cite{zhang2004configuration,biere2009handbook,khudabukhsh2016satenstein}.


However, research on {local search \PL{solvers} for general ILP is pretty rare. Some attempts on sub-classes of ILP have been proposed, such as over-constrained ILP \citep{walser1998domain,walser1998integer} and 0-1 ILP \citep{souza2014local, umetani2017exploiting, lei2021efficient, iser2023oracle}. \MC{ In \cite{prestwich2008constructive} authors proposed a local search algorithm that accepts general ILP form, but it is only tested on one special type of problem. There has been a work \cite{luteberget2023feasibility} that proposes an local search algorithm \textit{Feasibility Jump} for ILP that is tested on heterogeneous benchmark and won the 1st place in MIP 2022 Competition\footnote{https://www.mixedinteger.org/2022/competition/}. But \textit{Feasibility Jump} aims at finding feasible solutions and does not consider the objective function, which does not treat the original optimization purpose of ILP.} In general, the ability of local search for in quickly computing good solutions of ILP remains to be explored. 

\subsection{Contributions}

\PL{We develop a new local search algorithm Local-ILP, which is efficient for solving general ILP that is validated on a variety of different problems (i.e., the MIPLIB dataset).}

{Our algorithm is developed with two new operators tailored for ILP and a three-mode search framework. 
The two operators, namely the \textit{tight move} and \textit{lift move}, are designed based on a characterization of the solution space of ILP that involves the concept of boundary solutions. We show that searching for boundary solutions is complete for feasible and finite ILP instances, i.e., all optimal solutions belong to boundary solutions.}
{We also design a local search framework that switches between three modes, namely \search{}, \improve{}, and \restore{} modes. Depending on the state of the best-found solution and the current solution, the framework selects the process of the appropriate mode to execute, where appropriate operators are leveraged to improve the quality of the current solution according to different situations.}

For the \search{} and \restore{} modes, the \textit{tight move} operator is adopted, which jointly considers variables and constraints' information and tries to make some constraint tight. To distinguish important constraints and help guide the search to find feasible and high-quality solutions, we design a tailored weighting scheme and scoring function to select operations.
For the \improve{} mode, an efficient operator \textit{lift move} is used to improve the quality of the objective function while maintaining feasibility. To drive \textit{lift move} operation, we propose a way to compute a variable’s new candidate values called local domain reduction. Additionally, we also design a specialized scoring function for \textit{lift move} to pick a good operation.
By putting these together, we develop a local search ILP \PL{solver} called Local-ILP.

Experiments are conducted to evaluate Local-ILP on the MIPLIB benchmark, in which the ILP instances labeled hard and open are selected. 
We compare our algorithm with the state-of-the-art non-commercial ILP solver SCIP, as well as the state-of-the-art commercial solver Gurobi.
Experimental results show that, {in terms of finding a good feasible solution} within a reasonable time, Local-ILP is competitive and complementary with Gurobi, and significantly better than SCIP. 
{Experiments also  demonstrate that Local-ILP significantly outperforms another local search algorithm \textit{Feasibility Jump}, which won 1st place in MIP 2022 Workshop's Computational Competition.}
Moreover, Local-ILP establishes 6 new records for MIPLIB open instances by finding the new best solutions. 

Finally, we perform a theoretical analysis of our operators and algorithm based on the concept of boundary solutions, which are more explicit and closely related to the original form of ILP than the integral hull in polyhedral theory. We show that all feasible solutions visited by our algorithm are boundary solutions, efficiently avoiding unnecessary regions. 


\subsection{Outline}
The remainder of the paper is organized as follows: 
Section \ref{Preliminary} introduces the ILP problem and basic concepts for the local search algorithm, and then presents some basic characterizations of ILP we leveraged to design our local search algorithm. 
Section \ref{Our Local Search Framework for ILP} proposes our local search framework for solving general ILP. 
In sections \ref{Tight Move} and \ref{Lift Move}, we introduce two new \PL{operators} that are key techniques for different modes of the framework. 
Section \ref{Local-ILP Algorithm} presents detailed descriptions of how the Local-ILP algorithm implements the framework. 
The experimental results of our algorithm on public benchmark instances are reported in Section \ref{Experiments}. 
The analysis of our algorithm is presented in Section \ref{theory}, followed by the conclusion and future work in Section \ref{Conclusions}.

\section{Preliminary
\label{Preliminary}
}
In this section, we present some fundamental integer linear programming and local search concepts pertinent to the paper.

\subsection{Integer Linear Programming Problem}
An instance of generalized ILP has the following form:
\begin{equation} \label{eq1}
\begin{split}
Minimize\ \ \    &obj(\bl{x})=\bl{c}^\top \bl{x} \\ 
subject\ to\ \ \ \    &\bl{A}\bl{x} \leq \bl{b} \\
                       &\bl{x^l} \leq  \bl{x} \leq \bl{x^u}  \\
                       &\bl{x} \in \mathbb{Z}^n  \\
\end{split}
\end{equation}
over integer variables $\bl{x}$, where $\bl{A} \in \mathbb{R}^{m \times n}$, $\bl{b} \in \mathbb{R}^m$, $\bl{c} \in \mathbb{R}^n$, $\bl{x^l} \in (\mathbb{Z}\cup -\infty)^n$ and $\bl{x^u} \in (\mathbb{Z}\cup +\infty)^n$ are given inputs. The goal of ILP is to minimize the value of the objective function, while satisfying all constraints. We denote the \textit{i}-th constraint in the constraint system $\bl{A}\bl{x} \leq \bl{b}$ by $con_i$: $\bl{A_i}  \bl{x} \leq \bl{b_i}$, where $\bl{A_i} \in \mathbb{R}^{n}$, $\bl{b_i} \in \mathbb{R}$. The coefficient of $x_j$ in $con_i$ is $A_{ij}$ and $con_i$ contains $x_j$ if $A_{ij} \neq 0$. The variables' bounds are denoted by $\bl{x^l} \leq  \bl{x} \leq \bl{x^u}$; they are indeed parts of $\bl{A}\bl{x} \leq \bl{b}$, but will be treated separately from the coefficient matrix $\bl{A}$ in practical algorithms. The infinite value of $\bl{x^l}$ or $\bl{x^u}$ indicates no lower or upper bound on the corresponding variable. 

A solution $\bl{\alpha}$ for an ILP instance $Q$ is a vector of values assigned for each variable.
Specifically, $\alpha_j$ denotes the value of $x_j$. 
A solution $\bl{\alpha}$ of $Q$ satisfies $con_i$ in $\bl{A}\bl{x} \leq \bl{b}$ if $\bl{A_i} \cdot \bl{\alpha} \leq \bl{b_i}$, and $\bl{\alpha}$ is feasible if and only if it satisfies all constraints in $Q$. 
The value of the objective function of a solution $\bl{\alpha}$ is denoted as $obj(\bl{\alpha})$.

\subsection{Local Search Algorithm }
The local search algorithm is well used in solving combinatorial optimization problems. It explores the search space comprised of all solutions, each of which is a candidate solution. Normally, local search starts with a solution and iteratively alters the solution by modifying the value of one variable, in order to find a feasible solution with a high-quality objective function value. 
The key to local search is to decide, under different circumstances, which variables to modify and to what value, to get a new solution.

An \textbf{operator} defines how the candidate solution is modified, i.e., given variables to be modified, how to fix them to new values. 
When an operator is instantiated by specifying the variable to operate on, we obtain an \textbf{operation}.
Given an operation \op{}, the scoring function $score(op)$ is used to measure how good \op{} is. An operation \op{} is said to be {\bf positive} if $score(op) > 0$, which indicates that performing \op{} can improve the quality of the current solution.

\subsection{New Characterizations of ILP for Local Search} \label{convex}

Currently, the most common theoretical tool to analyze ILP is polyhedral theory \citep{schrijver1998theory}, of which the key component is the convex hull of all feasible solutions of an ILP (i.e., the integral hull), and it has established the equivalence between searching for the optimal solution and finding this integral hull. This equivalence {provides additional theoretical explanations} for the cutting plane method, which is widely used in commercial solvers, although the integral hull is hard to compute and could be done in exponential time.

However, the characterizations by integral hull are difficult to use for analyzing local search for ILP, as there is no explicit form of this integral hull and also no direct relations with constraints in the original form of ILP. \MC{To facilitate the analysis of local search algorithms for ILP, we introduce some new way to characterize solutions of ILP based on its original form, and show our local search operators and combinations are suited for ILP.}  
We briefly present here the intuitions that we leveraged to design our local search algorithm and will give precise formulations and theoretical arguments in Section~\ref{theory}.  

An important feature of ILP, is the linearity of its objective function and constraints. Although the ILP does not have true convexity, it still has some properties of this nature, such as all feasible solutions staying in the feasible domain of its LP relaxation, which is convex. 
From this aspect, some simple facts could be observed. Let $J=\{1, ..., n\}$, $\bl{e}_j$ be the unit vector with $1$ in $j$-th coordinate and $0$ in all other places. There is: 
\begin{fact}
    If for $\bl{x}_1, \bl{x}_2 \in \mathbb{Z}^n$, s.t. $\exists j \in J, k\in \mathbb{Z}, k> 0$, $\bl{x}_2 = \bl{x}_1+k\bl{e}_j$, $\bl{x}_1, \bl{x}_2$ are both feasible for an ILP instance, then for $\bl{x}'=\bl{x}_1+t\bl{e}_j$, $t\in \{0, ..., k\}$, $\bl{x}'$ is also feasible. 
\end{fact}
\begin{fact}
    If for $\bl{x}_1, \bl{x}_2 \in \mathbb{Z}^n$, s.t. $\exists j \in J, k\in \mathbb{Z}, k> 0$,  $\bl{x}_2 = \bl{x}_1+k\bl{e}_j$. Let $obj(\bl{x})$ be the objective function of an ILP. Then if $obj(\bl{x}_1)\leq obj(\bl{x}_2)$ then for $\bl{x}'=\bl{x}_1+t\bl{e}_j$, $t\in \{0, ..., k\}$, $obj(\bl{x}_1)\leq obj(\bl{x}') \leq obj(\bl{x}_2)$.  
\end{fact}

In other words, solutions that lie between two solutions that are different in only one dimension have the objective value that is also in between (or could be equal to) the two solutions. 

\MC{Thus we could make the following observations about ILP feasible solutions and optimal solutions, which will be used later to design new operators and our three-mode framework}: 

(1) all feasible solutions lie within a region (the integral hull) 

(2) all optimal solutions lie in the ``boundary'' of the region

Furthermore, we call a set of solutions a \textbf{complete} space to an ILP if all optimal solutions are contained in this space. We will formalize the concept of ``boundary'' in Section~\ref{theory} and show that boundary solutions are complete for an ILP. \MC{We design our operators based on this concept such that all feasible solutions our algorithm visits are boundary solutions.}

For a search algorithm, we can distinguish three different stages: 
(1) {enter} the integral hull (= find a feasible solution)

(2) find good solutions within the integral hull (= improve the quality of a feasible solution)

(3) get back to the integral hull when jumping out of it (= from an infeasible solution to get a good feasible solution)

\MC{The characters of the three stages motivate us to propose a three-mode framework, enabling our algorithm to have different behaviors in different situations, improving efficiency.}  

\MC{
In the following, we will first give descriptions of our algorithm, and later present the formal definitions and analysis of our algorithm in Section~\ref{theory}. 
}

\section{A New Local Search Framework for ILP
\label{Our Local Search Framework for ILP}
}

As depicted in Figure 1, we propose a new local search framework that takes advantage of adopting three different modes, namely \search{}, \improve{}, and \restore{} modes. In each mode, we use tailored operators to explore new solutions and thus diversify the search behaviors during different periods. 

\begin{figure}[ht]
\label{fig.1}
\centering
\centerline{\includegraphics[scale=0.16]{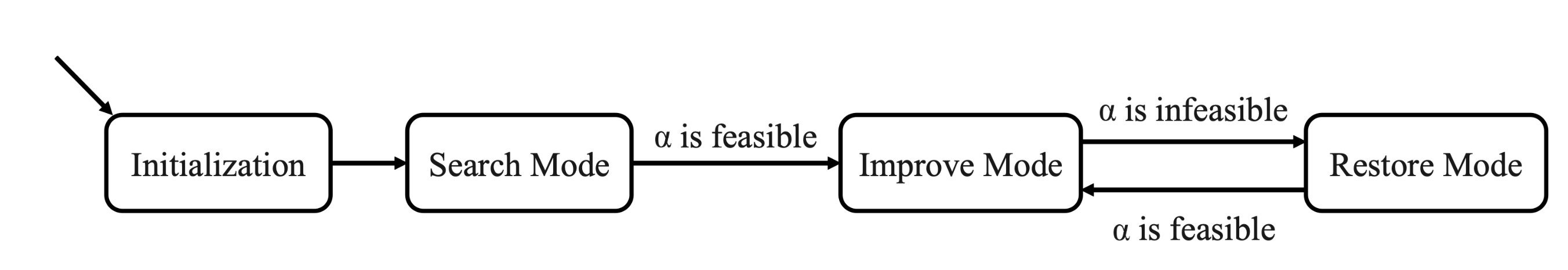}}
\caption{\centering A Local Search Framework for ILP}
\label{framework}
\end{figure}

(1) \textbf{Search mode}. The algorithm initially works in \search{} mode to find a feasible solution. When a feasible solution is found, the algorithm switches between \improve{} and \restore{} modes.

(2) \textbf{Improve mode}. Based on a known feasible solution, the \improve{} mode attempts to improve the value of the objective function while maintaining feasibility. If no such operation is found, the algorithm will break feasibility, obtain an infeasible solution, and enter \restore{} mode. 

(3) \textbf{Restore mode}. The \restore{} mode is to repair an infeasible solution to a good feasible solution. For an infeasible solution, the \restore{} mode concentrates on searching for a new high-quality feasible solution by considering the objective function and repairing the feasibility, and then returns to \improve{} mode after success.

\begin{algorithm}[htbp]
\caption{Local Search Framework for ILP}
\label{ls-1}
\KwIn{ILP instance $Q$, cut off time \textit{cutoff}}
\KwOut{A solution $\bl{\alpha}$ of $Q$ and its objective value}
$\bl{\alpha}$:= an initial solution; $\bl{\alpha^*}$:=$\emptyset$ ; \  $obj^*$:=$+\infty$;

\While{running\ time $<$ \textit{cutoff} } {
\lIf{$\bl{\alpha}$ is feasible \textbf{and} $obj(\bl{\alpha})<obj^*$}{
$\bl{\alpha^*}$:=$\alpha$; $obj^*$:=$obj(\bl{\alpha})$
}
\lIf{$\bl{\alpha^*}== \emptyset$}{
    perform operation for \search{} Mode 
}
\lElseIf {$\bl{\alpha}$ is feasible}{
        perform operation for \improve{} Mode
    }
\lElse{
perform operation for \restore{} Mode    
}
\lIf{enough steps to not improve}{
    restart the search}
}
\Return $\bl{\alpha^*}$  and $obj^*$ \;
\end{algorithm}

An outline of the framework is presented in Algorithm \ref{ls-1}. In the beginning, the algorithm constructs a solution $\bl{\alpha}$, and initializes the best-found solution $\bl{\alpha^*}$ as $\emptyset$ and its objective value to $+\infty$. 
The core of the algorithm consists of a loop (lines 2–7) in which solution $\bl{\alpha}$ is modified iteratively until a given time limit is reached. Depending on the state of the solutions $\bl{\alpha}$ and $\bl{\alpha^*}$, it chooses, \textit{Search}, \textit{Improve}, or \textit{Restore} modes to perform operations on variables (lines 4–6). Once a better feasible solution is discovered during the search, $\bl{\alpha^*}$ and $obj^*$ are updated correspondingly (line 3). The search is restarted if $\bl{\alpha^*}$ is not updated within a sufficient number of steps (line 7).
When the time limit is reached, the algorithm returns the best-found solution $\bl{\alpha^*}$, and its objective value $obj^*$. 

{Since we have a common general framework to integrate different operators in each mode, we describe the general framework with the notion of ``$X$ operation'' to represent the corresponding operation used in each mode.}
In each mode \textit{X} (\textit{X} is \search{}, \improve{} or \restore{}), an \textit{X} operation is iteratively picked to modify $\bl{\alpha}$, where an \textit{X} operation refers to an operation that is customized for mode \textit{X}. When the feasibility of $\bl{\alpha}$ is changed in the current mode, the algorithm shifts to another mode as we explained. 
All three modes adopt a general procedure as described in Algorithm \ref{ls-2}. It prefers to pick a positive operation (according to some heuristic) if one exists. If no such operation exists, at which point the algorithm is stuck in a local optimum, a random perturbation operation on $\bl{\alpha}$ is performed. 

\begin{algorithm}[ht]
\caption{Process for Mode $X$}
\label{ls-2}
\KwIn{ILP instance: $Q$, a solution: $\bl{\alpha}$}
\lIf{$\exists$ positive $X$ operations }{
    $op :=$  a positive $X$ operation 
}
\lElse{
     $op :=$  a random perturbation $X$ operation to escape local optima
}

perform \op{} to modify $\bl{\alpha}$;
\end{algorithm}

In the following two sections, we present those operations (also called \textit{X} operations) in each mode, which contain two new strategies raised by us, namely \textit{tight move} and \textit{lift move}, and section \ref{Local-ILP Algorithm} presents the whole Local-ILP algorithm that implements the above framework. 
\section{Tight Move Operator
\label{Tight Move}
}

In this section, we propose the \textit{tight move} (\textit{tm} for short) operator, which is used in \search{} mode and \restore{} mode, with an efficient weighting scheme and a tailored scoring function.

\subsection{Tight Move}

\PL{For solving ILP that may have infinite domains, a naive operator of local search is to modify the value of a variable $x_j$ by a fixed incremental value $inc$, i.e., $\alpha_j:=\alpha_j\pm inc$.
However, the fixed $inc$ has drawbacks that it is hard to choose and may be instance-dependent: if $inc$ is too small, it may take many iterations before making any violated constraints satisfied; if $inc$ is too large, the algorithm may even become so troublesome that it can never satisfy some constraints, and therefore be hard to find a feasible solution. Our \textit{tight move} operator modifies variable values according to the current solution and constraints and thus automatically adapts to different situations, }which is defined below.

\begin{myDef}\label{def:TM}
Given a solution $\bl{\alpha}$, the \textbf{tight move operator}, denoted as $tm(x_j, con_i, \bl{\alpha})$, assigns an integer variable $x_j$ to the value making the constraint $con_i$ as tight as possible and keeping the $x_j$'s bounds satisfied. Here, $con_i$ contains $x_j$ and could be either violated or satisfied. Precisely, let $\Delta=\bl{b_i}-\bl{A_i} \cdot \bl{\alpha}$, {which is commonly known as a slack}, a \tm{} operation is:

 (1) If $\Delta < 0$: $con_i$ is \textbf{violated}, there exists a {\it tm} operation $tm(x_j, con_i, \bl{\alpha})$ for variable $x_j$:

If $A_{ij} < 0$, then $tm(x_j, con_i, \bl{\alpha})$ increases $\alpha_j$ by $min(\left \lceil {\Delta}/{A_{ij}} \right \rceil ,  x^u_j-\alpha_j )$   

If $A_{ij} > 0$, then $tm(x_j, con_i, \bl{\alpha})$ decreases $\alpha_j$ by $min(\left| \left \lfloor {\Delta}/{A_{ij}}  \right \rfloor \right|, \left| x^l_j-\alpha_j \right|)$

 (2) If $\Delta \geq 0$: $con_i$ is \textbf{satisfied}, there exists a {\it tm} operation $tm(x_j, con_i, \bl{\alpha})$ for variable $x_j$:

 If $A_{ij} < 0$, then $tm(x_j, con_i, \bl{\alpha})$ decreases $\alpha_j$ by $min(\left| \left \lceil {\Delta}/{A_{ij}} \right \rceil \right| ,\left|  x^l_j-\alpha_j \right|)$
 
 If $A_{ij} > 0$, then $tm(x_j, con_i, \bl{\alpha})$ increases $\alpha_j$ by $min(\left \lfloor {\Delta}/{A_{ij}} \right \rfloor , x^u_j-\alpha_j )$

\end{myDef}

Our \textit{tight move} operator has two properties while keeping the  $x_j$’s bounds satisfied: 

(1) If $con_i$ is violated, it makes $con_i$ as close to being satisfied as possible while minimally influencing other constraints, as it selects the minimal change of $x_j$ to make $con_i$ do so. 

(2)If $con_i$ is satisfied, it can push $x_j$ to its extreme value, which keeps $con_i$ satisfied. This explores the influence of the maximal change of $x_j$ to the objective function and also helps to jump out of the local optimum. 

\PL{
The \textit{tight move} operator, which considers both the information of variables and constraints, takes the idea to modify a variable's value to make a constraint tight. 
This idea resembles some previous works which also consider some form of tightening constraints. The Simplex algorithm \citep{dantzig1955generalized} adopts similar operator and is commonly used to solve linear systems. 
Specifically, the \textit{tight move} is inspired by an operator named \textit{critical move}, which is proposed in the context of solving Satisfiability Modulo Integer Arithmetic Theories \citep{cai2022local}.
The \textit{critical move} operator assigns a variable $x_j$ to the threshold value making a violated literal (similar to a constraint in ILP) true, while our \textit{tight move} operator keeps variables satisfying their global bounds, and allows modifications from satisfied constraints, thus expanding the range of operations to choose from. 
}

\PL{
We will show in Section~\ref{theory} that although ILP does not have classical convexity, our \textit{tight move} operator still has good theoretical properties to produce promising solutions and avoid visiting unnecessary search space. 
}

\subsection{Scoring Function for Tight Move}
During the local search process, scoring functions are used to compare different operations and pick one to execute in each step.
For the \textit{tight move}, we propose a weighted scoring function.
Our scoring function has two ingredients: a weighting scheme to distinguish important constraints, and score computations to select operations. 

\subsubsection{Weighting Scheme}

In order to guide the search process, weighting techniques are widely applied in local search algorithms, and are used primarily to increase the weight of violated constraints, hence guiding the search process towards satisfying violated constraints. 

Here we present the weighting scheme we utilize, which will be further used in the scoring function for selecting operations. 

We assign an integral weight to each constraint and also the objective function, which are denoted as $w(con_i)$ and $w(obj)$, respectively.
At the beginning of the search, these weights are initialized to 1. To prevent too large values, we set an upper limit $ul_{con}$ and $ul_{obj}$ to them, respectively:
 
 (1)\PL{For the weight of the constraints, $ul_{con}=max(1000 , nCons)$, where $nCons$ denotes the number of constraints in the instance.}
 
 (2)\PL{For the weight of the objective function, $ul_{obj}=ul_{con}/10$.}


For weight setting, only when a solution is feasible do we consider it meaningful to improve the quality of the objective function. With this consideration, $w(obj)$ should not be too large compared to the weights of the constraints. 

\PL{Our weighting scheme is based on the probabilistic version of the PAWS scheme \citep{thornton2004additive, cai2013local, cai2022local}, which updates the weights according to a smoothing probability $sp$, and we used the same setting as \cite{cai2022local} to set $sp = 0.0003$.}
When the weighting scheme is activated, the weights of constraints and objective function are updated as follows:

(1) Once the weighting scheme is activated, the weights of the constraints are updated:
with probability $1 -sp$, for each violated constraint $con_i$, if $w(con_i) < ul_{con}$, $w(con_i)$ is increased by one; 
 otherwise, for each satisfied constraint $con_i$, if $w(con_i)>1$, $w(con_i)$ is decreased by one.

(2) The weight of the objective function is updated only if any feasible solution is found:
 with probability $1 -sp$, if $obj(\bl{\alpha}) \geq obj^*$ and $w(obj)<ul_{obj}$, $w(obj)$ is increased by one;
 otherwise, if $obj(\bl{\alpha}) < obj^*$ and $1 < w(obj)$, $w(obj)$ is decreased by one.

By changing the weights of constraints and thus focusing on constraints that are often violated in local optima, we help the local search process find feasible solutions. The weight updates for the objective function help guide the search towards solutions with better objective values.

\subsubsection{Score Computations}

Based on the weighting scheme, we propose the scoring function for \textit{tight move}, which helps the local search algorithm pick a \tm{} operation to execute. 
Specifically, it has two parts: score for reducing the violation of constraints, and score for improving the objective function. 

{\textbf{Score for reducing the violation of constraints}. The score for reducing the violation of constraints is determined according to the satisfiability of the constraints in the following three cases: 

(1) If a violated constraint $con_i$ is satisfied after performing \op{}, it incurs a positive reward of $w(con_i)$. 

(2) If a satisfied constraint $con_i$ is violated after performing \op{}, it incurs a negative reward of $-w(con_i)$. 

(3) If a violated constraint $con_i$ is still violated after performing \op{}, but $\bl{A_i} \cdot \bl{\alpha}$ is reduced, i.e., $con_i$ is closer to being satisfied, it incurs a positive reward of $\beta \cdot w(con_i)$; otherwise, if $\bl{A_i} \cdot \bl{\alpha} $ is increased and exacerbates $con_i$'s violation, the negative reward is $\beta \cdot -w(con_i)$, where $\beta \in \left[0,1 \right] $ is a parameter. 

The score of this impact of a \tm{} operation \op{}, denoted by $score_{reduce}(op)$, is set as the sum of the reward over all constraints obtained by performing \op{}.}

\textbf{Score for improving the objective function}. The score for improving the objective function of a \tm{} operation \op{} is denoted by $score_{improve}(op)$. If the value of the objective function after performing \op{} is smaller than $obj^*$, $score_{improve}(op):=w(obj)$; otherwise, $score_{improve}(op):=-w(obj)$. Note that $score_{improve}(op)$ is only meaningful in \restore{} mode, and is a constant in \search{} mode.
\begin{myDef}
    The \textbf{tight move score} of a \textit{tight move} operation \op{}, denoted as $score_{tm}(op)$, is defined as $$score_{tm}(op)=score_{reduce}(op)+score_{improve}(op)$$
\end{myDef}




\section{Lift Move Operator}
\label{Lift Move}
In this section, we propose the \textit{lift move} (\textit{lm} for short) operator, which is the key technique of the \improve{} mode. The property of the \textit{lift move} operator is that it can improve the quality of the objective function while maintaining feasibility. For this purpose, we propose a way to compute a variable's new candidate values called local domain reduction, and a specialized scoring function for \textit{lift move}.
\subsection{Local Domain Reduction}

To maintain the feasibility of a feasible solution, we must ensure that it satisfies every constraint. Therefore, we propose the local domain reduction to compute such a range to change a variable. 

\begin{myDef}
For a variable $x_j$ in a feasible solution $\bl{\alpha}$, its \textbf{local feasible domain}, denoted as $lfd(x_j,\bl{\alpha})$, is an interval for $x_j$ that when $x_j$ varies within this interval and all other variables stay unchanged, the satisfiability of all constraints will be kept unchanged. We call the local domain reduction the process to compute this local feasible domain $lfd(x_j,\bl{\alpha})$.
\end{myDef}

In order to compute $lfd(x_j,\bl{\alpha})$, we consider the feasible domain of $x_j$ in $con_i$, which is denoted as $ldc(x_j, con_i,\bl{\alpha})$ meaning the $x_j$ can vary within this interval while keeping the satisfiability of $con_i$, assuming other variables in $\bl{\alpha}$ keep unchanged, where $con_i$ is a constraint containing $x_j$. 
Specifically, let $\Delta=\bl{b_i}-\bl{A_i} \cdot \bl{\alpha}$, $ldc(x_j, con_i,\bl{\alpha})$  is derived according to the sign of $A_{ij}$: 

(1) If $A_{ij}\  < \ 0$, then $ldc(x_j, con_i,\bl{\alpha}) = \left[ \alpha_j + \left \lceil {\Delta}/{A_{ij}} \right \rceil , +\infty \right) $

(2) Otherwise,  $ldc(x_j, con_i,\bl{\alpha}) = \left(-\infty , \alpha_j +  \left \lfloor \Delta/{A_{ij}} \right \rfloor \right] $.

Then, with $ldc(x_j, con_i,\bl{\alpha})$, we calculate the local feasible domain of $x_j$ as follows: 
$$lfd(x_j,\bl{\alpha})=(\cap_i ldc(x_j, con_i,\bl{\alpha}) ) \cap [x_j^l, x_j^u]$$

\subsection{Lift Move}

Clearly, moving $x_j$ within integers of 
 $lfd(x_j,\bl{\alpha})$ does not break the feasibility of $\bl{\alpha}$ as long as the other variables are kept constant. So once the current solution $\bl{\alpha}$ is feasible, we can choose a reasonable integer in $lfd(x_j,\bl{\alpha})$ for updating $x_j$ to improve the quality of the objective function. We create the \textit{lift move} operator for this purpose, which is based on the local domain reduction. 

\begin{myDef}\label{def:LM}
For a feasible solution $\bl{\alpha}$, the \textbf{lift move operator}, denoted as $lm(x_j,\bl{\alpha})$, assigns $x_j$ to the upper or lower bound of $lfd(x_j,\bl{\alpha})$ to improve the objective function at most. Specifically, let $c_{j}$ denote the coefficient of $x_j$ in $\bl{c^\top} \bl{x}$, a \lm{} operation is described as follows: 

(1) If $c_j<0$, then $lm(x_j,\bl{\alpha})$ assigns $\alpha_j$ to the upper bound of $lfd(x_j,\bl{\alpha})$.

(2) If $c_j>0$, then $lm(x_j,\bl{\alpha})$ assigns $\alpha_j$ to the lower bound of $lfd(x_j,\bl{\alpha})$.

\end{myDef}

\subsection{Scoring Function for Lift Move}
If there are multiple variables in the objective function, then multiple \lm{} operations could be constructed. To guide the search in \improve{} mode, we customize a scoring function to select a \lm{} operation. Since all \lm{} operations will maintain the feasibility of the solution, we propose \textbf{lift score} to measure the improvement of the objective function.

\begin{myDef}
The \textbf{lift score} of a \textit{lift move} operation \op{}, denoted as $score_{lm}(op)$, is defined as  $$score_{lm}(op)= obj(\bl{\alpha})-obj(\bl{\alpha'})$$
where $\bl{\alpha}$ and $\bl{\alpha}'$ denotes the solution before and after performing \op{}.
\end{myDef}
In the \textit{Improve} mode, we pick the operation with the best $score_{lm}(op)$. 

\section{Full Local-ILP Algorithm with Three Modes
\label{Local-ILP Algorithm}
}

Based on the ideas in previous sections, we develop a local search \PL{ILP solver} called Local-ILP. As described in Section 3, after initialization, the algorithm works in three modes, namely \search{}, \improve{}, and \restore{} mode to iteratively modify $\bl{\alpha}$ until a given time limit is reached. This section is dedicated to the details of the initialization and the three modes of local search, as well as other optimization techniques.

{\bf Initialization}: Local-ILP generates a solution $\bl{\alpha}$, by assigning the variables one by one until all variables are assigned. As for a variable $x_j$, if $x^l_j>0$, it is assigned with $x^l_j$; if $x^u_j<0$, it is assigned with $x^u_j$. Otherwise, the variable is set to 0.

\subsection{Search Mode}

In \search{} mode (Algorithm \ref{ls-3}), the goal of the algorithm is to find the first feasible solution.

\begin{algorithm}[ht]
\caption{Process for \search{} Mode}
\label{ls-3}
\KwIn{ILP instance: $Q$, an infeasible solution: $\bl{\alpha}$}
\If{$\exists$ positive tm operation in violated constraints}{
    $op :=$  such an operation with the greatest $score_{tm}$ \;
}
\Else{
    update constraint weights \;
    $op :=$  a \tm{} operation with the greatest $score_{tm}$ in a random  violated constraint \;
}

perform \op{} to modify $\bl{\alpha}$;
\end{algorithm}

If the algorithm fails to find any positive \tm{} operations in violated constraints, it first updates the weights of constraints according to the weighting scheme described in Section 5 (line 6). Then, it picks a random violated constraint to choose a \tm{} operation with the greatest $score_{tm}$ (line 7). 

\subsection{Improve Mode}

The \improve{} mode (Algorithm \ref{ls-5}) seeks to improve the quality of the objective function's value while maintaining the feasibility of a feasible solution. If no such operation is found, the algorithm will break feasibility, obtain an infeasible solution, and enter \restore{} mode. 

\begin{algorithm}[ht]
\caption{Process for  Improve Mode}
\label{ls-5}
\KwIn{ILP instance: $Q$, a feasible solution: $\bl{\alpha}$ }
\If{$\exists$ positive lm operation
}{
    $op :=$  such an operation with the greatest $score_{lm}$ \;
} 
\lElse{
    $op :=$ a unit incremental move in the objective function within variables' global bounds
}
perform \op{} to modify $\bl{\alpha}$;
\end{algorithm}


If the algorithm fails to find any positive \lm{} operation, it randomly picks a variable $x_j$ in the objective function, and performs a simple \textbf{unit incremental move} in $x_j$ according to its coefficient $c_j$ (line 4-5), specifically, if $c_j<0$, then $\alpha_j=\alpha_j+1$; otherwise  $\alpha_j=\alpha_j-1$. If a unit incremental move will break the global bound of a variable, we randomly select another. A unit incremental move that keeps all global bounds must exist; otherwise, all variables are at the corresponding bound value that improves the objective function, which means the current solution is the optimal solution and we could finish the search.

\subsection{Restore Mode}

For an infeasible solution obtained from the \improve{} mode, the \restore{} mode (Algorithm \ref{ls-4}) focuses on repairing the feasibility to obtain a new high-quality feasible solution.

\begin{algorithm}[ht]
\caption{Process for Restore Mode}
\label{ls-4}
\KwIn{ILP instance: $Q$, an infeasible solution: $\bl{\alpha}$}
\If{$\exists$ positive tm operation in violated constraints}{
    $op :=$  such an operation with the greatest $score_{tm}$ \;
}
\ElseIf{$\exists$ positive tm operation in satisfied constraints}{
        $op :=$  such an operation with the greatest $score_{tm}$ \;
    
}
\Else{
    update constraint weights \;
    $op :=$  a \tm{} operation with the greatest $score_{tm}$ in a random  violated constraint \;
}

perform \op{} to modify $\bl{\alpha}$;
\end{algorithm}

If the algorithm fails to find any positive \tm{} operations in both violated and satisfied constraints, it updates weights and picks a random a \tm{} operation similarly to the \search{} mode (lines 6-7).

{The difference between \restore{} and \search{} modes is that \restore{} mode will pick \tm{} operations in satisfied constraints because the number of violated constraints in \restore{} mode is usually much less than in \search{} mode. More importantly, since $score_{improve}(op)$ is meaningful in \restore{} mode while a constant in \search{} mode, $score_{tm}$ used in \restore{} mode considers the quality of the objective function while \search{} mode does not.}

\subsection{Optimization Techniques}

The classical techniques of local search are applied in our algorithm, including random sampling, forbidding strategy, and restart mechanism.

{\bf Random Sampling}: We use a sampling method called Best from Multiple Selections, dubbed as BMS \citep{cai2015balance}. For all constraints that need to be considered in \search{} mode and \restore{} mode, the algorithm randomly samples a certain number of them, and for all operations derived from the sampled constraints that could be selected to apply, the algorithm randomly samples a certain number of these operations and selects the one with the highest $score_{tm}$. 
There are in total five parameters in the algorithm to control the number of samples, including:
$c_{v}$ for violated constraints, and $o_{v}$ for \tm{} operations in sampled violated constraints;
$c_{s}$ for satisfied constraints, and $o_{s}$ for \tm{} operations in sampled satisfied constraints;
$o_{r}$ for \tm{} operations in a random violated constraint. 


{\bf Forbidding Strategy}: 
We employ a forbidding strategy, the tabu strategy \citep{ElisaPappalardo2013HandbookOC,cai2022local}, to address the cycle phenomenon. After an operation is executed, the tabu strategy forbids the reverse operation in the following $tt$ iterations, where $tt$ is called tabu tenure, \PL{and we used the same setting as \cite{cai2022local} to set $tt = 3 + rand(10)$.} The tabu strategy is directly applied to Local-ILP. If a \tm{} operation that increases (decreases, resp.) the value of an integer variable $x_j$ is performed, then it is forbidden to decrease (increase, resp.) the value of $x_j$ by \tm{} operation in the following $tt$ iterations. 

{\bf Restart Mechanism:} \PL{The search is restarted when $\bl{\alpha^*}$ is not updated for enough steps iterations, which is simply set as 1500000.} At each restart, we crossover $\bl{\alpha^*}$ and random integers within variables’ bounds to reset  $\bl{\alpha}$: $x_j$ is assigned to $\alpha^*_j$  or a random integer within $[x^l_j, x^u_j]$ with 50\% probability of each, respectively. Additionally, all weights are restored to 1 when restarting.
\section{Experiments
\label{Experiments}
}

We carry out experiments to evaluate Local-ILP on the MIPLIB dataset\footnote{https://miplib.zib.de/}, which is a standard data set including a broad range of types of problems. 
{We compare our Local-ILP with state-of-the-art ILP solvers and another heuristic \textit{Feasibility Jump} in terms of their performance on finding a good feasible solution quickly.}
Also, experiments are conducted to analyze the effectiveness of the proposed new operators and framework. Additionally, the MIPLIB dataset records the best-known solutions for its instances, and we've established new records for 6 instances in the MIPLIB by Local-ILP. 

\subsection{Experiment Preliminaries}

\subsubsection{Implementation}
Local-ILP is implemented in C++ and compiled by g++ with the `-O3' option. There are two types of parameters: $\beta$ for the scoring function of \textit{tight move}, and the sampling numbers for the BMS heuristic. 
\PL{The default values of the parameters are shown in Table \ref{parameter_setting} tuned with 20\% randomly sampled instances.}

\begin{table*}[ht]
\setlength{\belowcaptionskip}{10pt}
\centering
\caption{\centering Tuned parameters of our proposed algorithms.
}
\label{parameter_setting}
\scalebox{0.8}
{
\begin{tabular}{c|c|c}
\hline
Parameter &Range &Final value \\
\hline
$\beta$ &\{0.1,0.2,...,1\} & 0.5    \\ 
$c_v$ &\{1,2,...,10\} & 3    \\                
$o_v$ &\{1000,2000,...,10000\} & 2000 \\ 
$c_s$ &\{10,20,...,100\} & 30 \\
$o_s$ &\{50,100,...,500\} & 350 \\
$o_r$ &\{50,100,...,500\} & 150 \\
\hline
\end{tabular}
}
\end{table*}

\subsubsection{Competitors} 
We compare Local-ILP with the latest state-of-the-art commercial and non-commercial ILP solvers, namely Gurobi 10.0.0\footnote{https://www.gurobi.com/} and SCIP 8.0.1\footnote{https://www.scipopt.org/}. For Gurobi, we use both its complete and heuristic versions. 
{Moreover, we compare Local-ILP with another heuristic \textit{Feasibility Jump}\footnote{https://github.com/sintef/feasibilityjump} (FJ for short), which won 1st place in MIP 2022 Workshop’s Computational Competition.}
For all of the competitors, we always use their default parameter settings.

\subsubsection{Benchmarks}\label{bench}
MIPLIB is widely recognized as the standard benchmark set for integer linear programming problems.
Our experiments are carried out with the union of MIPLIB 2003 \citep{achterberg2006miplib}, MIPLIB 2010 \citep{koch2011miplib} and MIPLIB 2017 \citep{gleixner2021miplib}, selecting the ILP instances tagged hard and open.
As Local-ILP is not an exhaustive algorithm, infeasible instances are excluded, resulting in a benchmark consisting of 121 instances.

\subsubsection{Experiment Setup}
All experiments are carried out on a server with AMD EPYC 7763 CPU and 2048G RAM under the system Ubuntu 20.04.4. For each instance, each solver is executed by one thread with time limits of 10, 60, and 300 seconds. 

we use 3 metrics to evaluate the performance of each solver:

    (1) \textbf{\feas{}}. To evaluate the ability to find a feasible solution for an instance, we count \feas{}, the number of instances where a solver can find a feasible solution within the time limit.
    
    (2) \textbf{\win{}}. To evaluate the ability to find a high-quality feasible for an instance, we count \win{}, the number of instances where a solver finds the best solution among all solutions output by tested solvers within the time limit. Note that, for each instance, if the best solution is found by more than one solver, this instance is counted as \win{} for all these solvers; and if all solvers find no solution, this instance is not counted as \win{} for any solver. 
    
    (3) {\textbf{\textit{P(T)}}. We also use a well-established measure, the primal integral \textit{P(T)} \citep{berthold2013measuring}, to evaluate the performance of each solver, which depends on the quality of solutions found during the solving process as well as on the points in time when they are found. To define the primal integral \textit{P(T)}, we first define the primal gap $\gamma$ and the primal gap function $p(t)$.}
    
    {Let $\bl{\tilde{x}}$ denote a solution for an ILP instance, and $\bl{\tilde{x}}_{opt}$ denote an optimal (or best known) solution for that instance, primal gap $\gamma$ is defined as \begin{equation}
        \gamma(\bl{\tilde{x}})=\left\{
         \begin{array}{ll}
         0, & if \  |\bl{c^\top} \bl{\tilde{x}}_{opt}| = |\bl{c^\top} \bl{\tilde{x}}| = 0,\\
         1, & if \ \bl{c^\top} \bl{\tilde{x}}_{opt} \cdot \bl{c^\top} \bl{\tilde{x}} < 0,\\
         \frac{|\bl{c^\top} \bl{\tilde{x}}_{opt} - \bl{c^\top} \bl{\tilde{x}}|}{max\{|\bl{c^\top} \bl{\tilde{x}}_{opt}|, |\bl{c^\top} \bl{\tilde{x}}|\}},& else.  
         \end{array}
        \right.
    \end{equation}}

    {Let $t_{max} \in \mathbb{R}_{\geq 0}$ be the time limit of an ILP solver. Its primal gap function $p(t)$: $[0, t_{max}] \rightarrow [0, 1]$, is defined as: \begin{equation}
        p(t)=\left\{
         \begin{array}{ll}
         1, & if \ no\  feasible\ solution\  until\  point\  t,\\
         \gamma(\bl{\tilde{x}}(t)), &\bl{\tilde{x}}(t)~is~the\  best\ found\  solution\  at\  point\ t,\  else.\\
         \end{array}
        \right.
    \end{equation}}
    
    {Let $T \in [0, t_{max}]$ and let $t_i \in [0, T]\ for\ i \in 1,..., I-1$ be the points in time when a new best solution is found, $t_0 = 0, t_I = T$. We define the primal integral $P(T)$ of a run as:
    \begin{equation}
        P(T):= \int_{t=0}^{T}p(t)dt =\sum_{i=1}^{I}p(t_{i-1})\cdot (t_i-t_{i-1})
    \end{equation}}

{For a benchmark, the primal integral $P(T)$ is the average primal integral of each instance. According to the above definition, smaller values of $P(T)$ indicate better performance.}

We organize the overall results into three types (\win{}, \feas{}, \pt{}). 
For each time limit setting, the best results are marked in \textbf{bold}.

\subsection{Comparisons with Competitors}\label{8.2}

In the MIPLIB dataset, each instance may contain multiple types of constraints\footnote{https://miplib.zib.de/statistics.html}, such as knapsack constraints, set covering constraints, and so on. We categorize all instances by the type of their main constraint class. 
We mark the instance as hybrid if it contains multiple main constraint classes.

The results of the comparison with all competitors in terms of the ability to find a feasible solution, the quality of the best-found solution, and the primal integral are shown in Tables \ref{com-res-feas}, \ref{com-res-win} and \ref{com-res-primalIntegral-1}, respectively.

\begin{table*}[ht]
\centering
\setlength{\belowcaptionskip}{10pt}
\caption{\centering Empirical results on comparing Local-ILP with FJ, SCIP, and Gurobi in terms of the ability to find a feasible solution, with 10s, 60s, and 300s time limits. \#inst denotes the number of instances in each class.
}
\label{com-res-feas}
\renewcommand\arraystretch{1.5}
\scalebox{0.45}
{
\begin{tabular}{l|c|ccccc|ccccc|ccccc}
\hline
\multirow{4}{*}{ Main Constraint Class}&\multirow{4}{*}{\#inst}& \multicolumn{12}{c}{\feas{}}\\ \cline{3-17} 
&  & \multicolumn{5}{c|}{10 s}& \multicolumn{5}{c|}{60 s}& \multicolumn{5}{c}{300 s}\\ \cline{3-17}
& & Local-ILP&FJ& SCIP& \multicolumn{2}{c|}{Gurobi}& Local-ILP&FJ& SCIP      & \multicolumn{2}{c|}{Gurobi}& Local-ILP&FJ& SCIP & \multicolumn{2}{c}{Gurobi}\\ 
\cmidrule(lr){6-7} \cmidrule(lr){11-12} \cmidrule(lr){16-17} 
&&&&& comp.& heur.&&&& comp.& heur.&&&& comp.& heur.\\ \hline
Singleton &2 &\textbf{2} &0 &\textbf{2} &\textbf{2} &\textbf{2} &\textbf{2} &0 &\textbf{2} &\textbf{2} &\textbf{2} &\textbf{2} &0 &\textbf{2} &\textbf{2} &\textbf{2} \\
Aggregations &2 &\textbf{1} &\textbf{1} &0 &\textbf{1} &\textbf{1} &\textbf{1} &\textbf{1} &\textbf{1} &\textbf{1} &\textbf{1} &\textbf{1} &\textbf{1} &\textbf{1} &\textbf{1} &\textbf{1} \\
Bin Packing &2 &1 &1 &0 &\textbf{2} &\textbf{2} &\textbf{2} &\textbf{2} &1 &\textbf{2} &\textbf{2} &\textbf{2} &\textbf{2} &1 &\textbf{2} &\textbf{2} \\
Equation Knapsack &3 &0 &0 &0 &0 &0 &0 &0 &0 &0 &0 &0 &0 &0 &0 &0 \\
Knapsack &4 &\textbf{4} &\textbf{4} &3 &3 &3 &\textbf{4} &\textbf{4} &3 &3 &3 &\textbf{4} &\textbf{4} &3 &\textbf{4} &\textbf{4} \\
Set Packing &5 &\textbf{4} &\textbf{4} &3 &\textbf{4} &\textbf{4} &\textbf{4} &\textbf{4} &\textbf{4} &\textbf{4} &\textbf{4} &\textbf{4} &\textbf{4} &\textbf{4} &\textbf{4} &\textbf{4} \\
Cardinality &6 &\textbf{2} &1 &1 &\textbf{2} &\textbf{2} &\textbf{2} &\textbf{2} &1 &\textbf{2} &\textbf{2} &2 &2 &2 &\textbf{3} &\textbf{3} \\
Hybrid &7 &4 &\textbf{5} &2 &3 &3 &4 &\textbf{5} &3 &4 &4 &4 &\textbf{5} &4 &\textbf{5} &\textbf{5} \\
Mixed Binary &8 &\textbf{4} &\textbf{4} &1 &2 &2 &4 &\textbf{5} &1 &3 &3 &\textbf{5} &\textbf{5} &2 &3 &3 \\
Set Partitioning &9 &5 &5 &3 &\textbf{7} &\textbf{7} &6 &\textbf{7} &4 &\textbf{7} &\textbf{7} &\textbf{7} &\textbf{7} &5 &\textbf{7} &\textbf{7} \\
Set Covering &11 &\textbf{9} &\textbf{9} &7 &\textbf{9} &\textbf{9} &\textbf{9} &\textbf{9} &8 &\textbf{9} &\textbf{9} &\textbf{9} &\textbf{9} &8 &\textbf{9} &\textbf{9} \\
Precedence &13 &\textbf{12} &9 &11 &\textbf{12} &\textbf{12} &\textbf{12} &10 &\textbf{12} &\textbf{12} &\textbf{12} &\textbf{12} &10 &\textbf{12} &\textbf{12} &\textbf{12} \\
General Linear &15 &\textbf{10} &6 &6 &8 &8 &\textbf{10} &6 &7 &9 &\textbf{10} &10 &6 &9 &\textbf{11} &\textbf{11} \\
Variable Bound &16 &\textbf{14} &12 &10 &13 &13 &\textbf{15} &13 &12 &14 &14 &\textbf{15} &13 &12 &14 &14 \\
Invariant Knapsack &18 &\textbf{13} &12 &10 &12 &12 &\textbf{13} &12 &11 &12 &12 &13 &12 &13 &\textbf{15} &\textbf{15} \\
\hline
Total &121 &\textbf{85} &73 &59 &80 &80 &\textbf{88} &80 &70 &84 &85 &90 &80 &78 &\textbf{92} &\textbf{92} \\
\hline
\end{tabular}
}
\end{table*}

\subsubsection{The Ability to Find a Feasible Solution} 
Here we present the results of the number of instances in which a feasible solution could be found by each solver.  

As shown in Table~\ref{com-res-feas}, Local-ILP performs best with 12 types of all main constraint classes in the 10s time limit, 13 types in the 60s, and 10 types in the 300s. 
Local-ILP wins the most types in the 10s and 60s time limits, and the second most in the 300s. 

In general, for the ability to find a feasible solution quickly, SCIP is the worst solver with the least \feas{}. It is obvious that Local-ILP, FJ, and Gurobi perform better than SCIP. Local-ILP performs best in the 10s and 60s time limits but loses to Gurobi in the 300s. 

\subsubsection{ The Quality of the Best Found Solution } 
For each solver, we present the number of instances where it found the best solution among all candidate solvers. 
\begin{table*}[ht]
\centering
\setlength{\belowcaptionskip}{10pt}
\caption{\centering Empirical results on comparing Local-ILP with FJ, SCIP, and Gurobi in terms of the quality of the best-found solution, with 10s, 60s, and 300s time limits. \#inst denotes the number of instances in each class.
}
\label{com-res-win}
\renewcommand\arraystretch{1.5}
\scalebox{0.45}
{
\begin{tabular}{l|c|ccccc|ccccc|ccccc}
\hline
\multirow{4}{*}{ Main Constraint Class}&\multirow{4}{*}{\#inst}& \multicolumn{12}{c}{\win{}}\\ \cline{3-17} 
&  & \multicolumn{5}{c|}{10 s}& \multicolumn{5}{c|}{60 s}& \multicolumn{5}{c}{300 s}\\ \cline{3-17}
& & Local-ILP&FJ& SCIP& \multicolumn{2}{c|}{Gurobi}& Local-ILP&FJ& SCIP      & \multicolumn{2}{c|}{Gurobi}& Local-ILP&FJ& SCIP & \multicolumn{2}{c}{Gurobi}\\ 
\cmidrule(lr){6-7} \cmidrule(lr){11-12} \cmidrule(lr){16-17} 
&&&&& comp.& heur.&&&& comp.& heur.&&&& comp.& heur.\\ \hline
Singleton &2 &\textbf{1} &0 &0 &0 &\textbf{1} &\textbf{1} &0 &0 &0 &\textbf{1} &\textbf{1} &0 &0 &\textbf{1} &0 \\
Aggregations &2 &0 &0 &0 &0 &\textbf{1} &0 &0 &0 &0 &\textbf{1} &0 &0 &0 &\textbf{1} &0 \\
Bin Packing &2 &0 &0 &0 &\textbf{2} &\textbf{2} &0 &0 &0 &\textbf{2} &\textbf{2} &0 &0 &\textbf{1} &\textbf{1} &\textbf{1} \\
Equation Knapsack &3 &0 &0 &0 &0 &0 &0 &0 &0 &0 &0 &0 &0 &0 &0 &0 \\
Knapsack &4 &\textbf{2} &0 &1 &0 &1 &2 &0 &0 &0 &\textbf{3} &1 &0 &0 &\textbf{3} &2 \\
Set Packing &5 &1 &0 &0 &\textbf{3} &2 &\textbf{2} &0 &0 &1 &\textbf{2} &\textbf{2} &0 &0 &1 &\textbf{2} \\
Cardinality &6 &\textbf{1} &0 &0 &\textbf{1} &0 &\textbf{1} &0 &0 &0 &\textbf{1} &0 &0 &0 &1 &\textbf{3} \\
Hybrid &7 &\textbf{3} &2 &0 &0 &0 &\textbf{3} &1 &0 &1 &1 &2 &1 &0 &2 &\textbf{3} \\
Mixed Binary &8 &\textbf{2} &\textbf{2} &1 &\textbf{2} &\textbf{2} &2 &\textbf{3} &1 &1 &2 &\textbf{2} &1 &1 &0 &1 \\
Set Partitioning &9 &3 &0 &0 &3 &\textbf{4} &\textbf{3} &1 &1 &1 &2 &1 &1 &0 &3 &\textbf{4} \\
Set Covering &11 &\textbf{5} &0 &1 &4 &4 &2 &0 &2 &2 &\textbf{3} &1 &0 &0 &1 &\textbf{8} \\
Precedence &13 &\textbf{6} &0 &0 &3 &5 &\textbf{8} &0 &1 &2 &1 &\textbf{5} &0 &1 &1 &\textbf{5} \\
General Linear &15 &3 &0 &1 &5 &\textbf{6} &1 &0 &1 &5 &\textbf{7} &1 &0 &1 &\textbf{7} &\textbf{7} \\
Variable Bound &16 &\textbf{11} &1 &0 &1 &3 &\textbf{9} &1 &2 &3 &7 &6 &0 &1 &4 &\textbf{9} \\
Invariant Knapsack &18 &\textbf{10} &1 &2 &5 &5 &\textbf{8} &1 &1 &7 &6 &4 &0 &2 &\textbf{10} &8 \\
\hline
Total &121 &\textbf{48} &6 &6 &29 &36 &\textbf{42} &7 &9 &25 &39 &26 &3 &7 &36 &\textbf{53} \\
\hline
\end{tabular}
}
\end{table*}

As shown in Table~\ref{com-res-win}, Local-ILP performs best for 9 types of all 15 main constraint classes in the 10s and 60s time limit, and 4 types in the 300s time limit. 
Local-ILP wins the most types in the 10s and 60s time limits, and the second most in the 300s time limit. 
In particular, Local-ILP exhibits the best performance for the types of \textit{singleton}, \textit{mixed binary}, and \textit{precedence} main constraint classes over all time limits.

Overall, for the ability to find a high-quality solution quickly, FJ performs the worst with the least \win{} in each time limit setting. Local-ILP and Gurobi perform much better than FJ and SCIP. For the comparison between Local-ILP and Gurobi, Local-ILP consistently performs best in the 10s and 60s time limits, but in the 300s time limits, Gurobi wins more instances than Local-ILP, especially its heuristic version.
\begin{table*}[!ht]
\centering
\setlength{\belowcaptionskip}{10pt}
\caption{\centering Empirical results on comparing Local-ILP with FJ, SCIP, and Gurobi in terms of the primal integral $P(T)$, with 10s, 60s, and 300s time limits. \#inst denotes the number of instances in each class.
}
\label{com-res-primalIntegral-1}
\renewcommand\arraystretch{1.5}
\scalebox{0.45}
{
\begin{tabular}{l|c|ccccc|ccccc|ccccc}
\hline
\multirow{4}{*}{ Main Constraint Class}&\multirow{4}{*}{\#inst}& \multicolumn{12}{c}{\pt{}}\\ \cline{3-17} 
&  & \multicolumn{5}{c|}{10 s}& \multicolumn{5}{c|}{60 s}& \multicolumn{5}{c}{300 s}\\ \cline{3-17}
& & Local-ILP&FJ& SCIP& \multicolumn{2}{c|}{Gurobi}& Local-ILP&FJ& SCIP      & \multicolumn{2}{c|}{Gurobi}& Local-ILP&FJ& SCIP & \multicolumn{2}{c}{Gurobi}\\ 
\cmidrule(lr){6-7} \cmidrule(lr){11-12} \cmidrule(lr){16-17} 
&&&&& comp.& heur.&&&& comp.& heur.&&&& comp.& heur.\\ \hline
Singleton &2 &\textbf{0.274} &1.000 &0.423 &0.352 &0.373 &\textbf{0.176} &1.000 &0.370 &0.189 &0.196 &0.143 &1.000 &0.225 &0.135 &\textbf{0.126} \\
Aggregations &2 &0.802 &0.854 &1.000 &\textbf{0.538} &0.542 &0.699 &0.809 &0.772 &0.509 &\textbf{0.508} &0.632 &0.802 &0.557 &0.502 &\textbf{0.502} \\
Bin Packing &2 &0.989 &0.997 &1.000 &0.960 &\textbf{0.959} &0.964 &0.994 &0.999 &\textbf{0.752} &0.760 &0.954 &0.994 &0.987 &\textbf{0.604} &0.606 \\
Equation Knapsack &3 &1.000 &1.000 &1.000 &1.000 &1.000 &1.000 &1.000 &1.000 &1.000 &1.000 &1.000 &1.000 &1.000 &1.000 &1.000 \\
Knapsack &4 &\textbf{0.166} &0.431 &0.360 &0.320 &0.320 &\textbf{0.143} &0.427 &0.309 &0.281 &0.280 &\textbf{0.133} &0.427 &0.298 &0.216 &0.209 \\
Set Packing &5 &\textbf{0.367} &0.903 &0.715 &0.424 &0.423 &0.293 &0.888 &0.468 &0.275 &\textbf{0.274} &0.275 &0.885 &0.314 &0.233 &\textbf{0.233} \\
Cardinality &6 &\textbf{0.690} &1.000 &0.861 &0.751 &0.752 &\textbf{0.657} &0.965 &0.825 &0.721 &0.720 &0.651 &0.871 &0.785 &0.564 &\textbf{0.557} \\
Hybrid &7 &\textbf{0.641} &0.682 &1.000 &0.780 &0.780 &\textbf{0.540} &0.605 &1.000 &0.648 &0.648 &\textbf{0.524} &0.592 &0.671 &0.544 &0.542 \\
Mixed Binary &8 &0.894 &\textbf{0.886} &1.000 &0.986 &0.986 &0.874 &\textbf{0.874} &1.000 &0.976 &0.976 &\textbf{0.867} &0.871 &0.992 &0.973 &0.970 \\
Set Partitioning &9 &0.906 &0.901 &0.927 &\textbf{0.819} &0.819 &0.862 &0.866 &0.895 &\textbf{0.811} &0.812 &0.836 &0.855 &0.822 &0.752 &\textbf{0.744} \\
Set Covering &11 &0.734 &0.777 &0.757 &\textbf{0.623} &0.623 &0.727 &0.773 &0.679 &0.496 &\textbf{0.491} &0.723 &0.773 &0.549 &0.376 &\textbf{0.346} \\
Precedence &13 &\textbf{0.694} &0.993 &0.877 &0.729 &0.722 &\textbf{0.547} &0.993 &0.727 &0.625 &0.622 &\textbf{0.439} &0.993 &0.629 &0.513 &0.517 \\
General Linear &15 &0.516 &0.751 &0.696 &\textbf{0.480} &0.483 &0.503 &0.736 &0.607 &0.408 &\textbf{0.404} &0.501 &0.734 &0.508 &\textbf{0.301} &0.321 \\
Variable Bound &16 &\textbf{0.506} &0.982 &0.831 &0.624 &0.619 &\textbf{0.371} &0.981 &0.633 &0.479 &0.472 &0.335 &0.981 &0.464 &0.356 &\textbf{0.321} \\
Invariant Knapsack &18 &\textbf{0.670} &0.860 &0.754 &0.720 &0.720 &0.621 &0.858 &0.732 &\textbf{0.619} &0.628 &0.563 &0.857 &0.693 &\textbf{0.483} &0.496 \\
\hline
Total &121 &\textbf{0.658} &0.871 &0.807 &0.675 &0.675 &0.597 &0.860 &0.721 &0.588 &\textbf{0.587} &0.563 &0.854 &0.628 &0.492 &\textbf{0.487} \\
\hline
\end{tabular}
}
\end{table*}

\subsubsection{Primal Integral} 
{Here we present the results of the average primal integral $P(T)$ of each solver.}
{As shown in Table~\ref{com-res-primalIntegral-1}, Local-ILP performs best with 8 types of all main constraint classes in the 10s time limits, 6 types in the 60s, and 4 types in the 300s. 
Moreover, Local-ILP performs best in most types in the 10s and 60s time limits. 
Particularly, Local-ILP exhibits the best performance for the types of \textit{knapsack}, \textit{mixed binary}, \textit{hybrid}, and \textit{precedence} main constraint classes over all time limits.}

{In general, for the primal integral, Local-ILP performs best in the 10s time limit, and is competitive with Gurobi in the 60s with a very close $P(T)$, but obviously loses to it in the 300s. However, Local-ILP is significantly better than FJ and SCIP over all time limits.}

\subsection{Effectiveness of Proposed New Techniques} \label{8.4}
To analyze the effectiveness of our proposed new techniques in Local-ILP, we tested 6 variations of our Local-ILP algorithm as follows: 

(1) To analyze the effectiveness of the \textit{tight move} operator. We modify Local-ILP by replacing the \textit{tight move} operator with the operator that directly modifies an integer variable by a fixed increment \textit{inc}, leading to two versions v\underline{~}fix\underline{~}1 and v\underline{~}fix\underline{~}5, where \textit{inc} is set as 1 and 5, respectively.

(2) To analyze the effectiveness of the \textit{lift move} operator and \textit{Improve} mode, we modify Local-ILP by removing the \textit{Improve} mode from the framework and using only \textit{Search} and \textit{Restore} modes, leading to the version v\underline{~}no\underline{~}improve. 
    To analyze the effectiveness of the \textit{Restore} mode, we modify Local-ILP by removing the \textit{Restore} mode from the framework and returning to \search{} mode when \textit{Improve} mode stuck in a local optimum, leading to the version v\underline{~}no\underline{~}restore.

(3) To compare different ways to escape from the local optimum in \textit{Improve} mode. We modify Local-ILP by replacing the unit incremental move with operators that have larger step sizes, leading to two versions v\underline{~}per\underline{~}bound and v\underline{~}per\underline{~}random, where the step size is set as the distance to the bound of the variable and a random size between 1 and the bound, respectively.

\begin{table*}[ht]
\centering
\setlength{\belowcaptionskip}{10pt}
\caption{\centering Empirical results on comparing Local-ILP with 6 variations,  with 10s, 60s and 300s time limits.
}
\label{com-res-ablation}
\renewcommand\arraystretch{1.6}
\scalebox{0.5}
{
\begin{tabular}{c|cccccc|cccccc|cccccc}
\hline
\multirow{2}{*}{TimeLimit} 
& \multicolumn{2}{c}{Local-ILP}
& \multicolumn{2}{c}{v\underline{~}fix\underline{~}1}
& \multicolumn{2}{c|}{v\underline{~}fix\underline{~}5}
& \multicolumn{2}{c}{Local-ILP}
& \multicolumn{2}{c}{v\underline{~}no\underline{~}improve}
& \multicolumn{2}{c|}{ v\underline{~}no\underline{~}restore}
& \multicolumn{2}{c}{Local-ILP}
& \multicolumn{2}{c}{v\underline{~}per\underline{~}bound}
& \multicolumn{2}{c}{v\underline{~}per\underline{~}random} \\ 
 \cmidrule(lr){2-3}\cmidrule(lr){4-5}\cmidrule(lr){6-7}\cmidrule(lr){8-9}\cmidrule(lr){10-11}\cmidrule(lr){12-13}\cmidrule(lr){14-15}\cmidrule(lr){16-17}\cmidrule(lr){18-19}
 &\win{} &\pt{}&\win{} &\pt{}&\win{} &\pt{}&\win{} &\pt{}&\win{} &\pt{}&\win{} &\pt{}&\win{} &\pt{}&\win{} &\pt{}&\win{} &\pt{}\\ 
\hline
10 s  &\textbf{64} &\textbf{0.658} &27 &0.723 &5 &0.924 &\textbf{61} &\textbf{0.658} &7 &0.811 &37 &0.662 &\textbf{56} &\textbf{0.658} &27 &0.682 &37 &0.680    \\
60 s  &\textbf{71} &\textbf{0.597} &31 &0.675 &1 &0.910 &\textbf{63} &\textbf{0.597} &4 &0.797 &35 &0.608 &\textbf{54} &\textbf{0.597} &23 &0.627 &45 &0.618    \\
300 s &\textbf{70} &\textbf{0.563} &31 &0.649 &2 &0.902 &\textbf{67} &\textbf{0.563} &4 &0.793 &31 &0.585 &\textbf{58} &\textbf{0.563} &26 &0.602 &41 &0.576    \\
\hline
\end{tabular}
}
\end{table*}

We compare Local-ILP with these modified versions on the benchmark, with 10s, 60s, and 300s time limits.
As shown in Table~\ref{com-res-ablation}, Local-ILP outperforms all other variations, confirming the effectiveness of the strategies.

\subsection{Stability with Repetitive Experiments}
{To examine the stability of Local-ILP which involves randomness, we run Local-ILP 10 times with 10 different seeds on the benchmark for 10s, 60s, and 300s time limits.}

\begin{table*}[ht]
\centering
\setlength{\belowcaptionskip}{10pt}
\caption{\centering Experimental results of Local-ILP with 10 different seeds on the benchmark for 10s, 60s, and 300s time limits.}
\label{com-res-stable}
\scalebox{0.65}{
    \begin{tabular}{ccc|ccc|ccc}
    \hline
     \multicolumn{3}{c|}{10 s}
    & \multicolumn{3}{c|}{60 s}
    & \multicolumn{3}{c}{300 s}
    \\
     \cmidrule(lr){1-3}  \cmidrule(lr){4-6} \cmidrule(lr){7-9}
     $avg_{P(T)}$ & $std_{P(T)}$ & ${std_{P(T)}}/{avg_{P(T)}}$ &
     $avg_{P(T)}$ & $std_{P(T)}$ & ${std_{P(T)}}/{avg_{P(T)}}$ &
     $avg_{P(T)}$ & $std_{P(T)}$ & ${std_{P(T)}}/{avg_{P(T)}}$ \\
    \hline
    0.65634& 0.00278& 0.00424 &0.59599& 0.00462& 0.00775 &0.56427& 0.00592& 0.01049\\
    
    \hline
    \end{tabular}
}
\end{table*}

{For all 10 times, we denote the average primal integral of each time by $avg_{P(T)}$, and the standard deviation of the primal integral by $std_{P(T)}$. The experimental results presented in Table~\ref{com-res-stable} demonstrate that, for each time limit, the values of ${std_{P(T)}}/{avg_{P(T)}}$ for Local-ILP are less than 1.1\%, indicating Local-ILP exhibits stable performance.}

\subsection{New Records to Open Instances}

The optimal solutions for the instances labeled open have not yet been reported. The current best solutions known for each open instance are available on the official website of MIPLIB 2017. 

\begin{table*}[ht]
\centering
\setlength{\belowcaptionskip}{10pt}
\caption{\centering New records to open instances
}
\label{record}
\renewcommand\arraystretch{1.5}
\scalebox{0.65}
{
\begin{tabular}{l|lcl}
\hline
Instance      & Local-ILP          & \multicolumn{1}{l}{Previous Objective} & Constraint Classification                                \\ \hline
sorrell7      & \textbf{-197}  & -196                                   & variable bound                                           \\
supportcase22 & \textbf{117}  & N/A                                    & set covering, aggregations, bin packing and mixed binary \\
cdc7-4-3-2    & \textbf{-294}  & -289                                   & set packing                                              \\
ns1828997     & \textbf{8}    & 9                                      & precedence, invariant knapsack, variable bound and cardinality                                  \\
scpm1         & \textbf{544}  & 554                                    & set covering                                             \\
scpn2         & \textbf{490}  & 501                                    & set covering                                             \\ \hline
\end{tabular}
}
\end{table*}

As shown in Table \ref{record}, Local-ILP established the new best-known objective values for 6 instances. 
These 6 instances contain different constraint types, which simultaneously demonstrate the strong solving power of Local-ILP and its applicability to diverse types of problems. 
{The new records have been submitted to MIPLIB 2017 and have been accepted; the links to the website are denoted in the footnotes
\footnote{https://miplib.zib.de/instance\underline{~}details\underline{~}sorrell7.html}
\footnote{https://miplib.zib.de/instance\underline{~}details\underline{~}supportcase22.html}
\footnote{https://miplib.zib.de/instance\underline{~}details\underline{~}cdc7-4-3-2.html}
\footnote{https://miplib.zib.de/instance\underline{~}details\underline{~}ns1828997.html} \footnote{https://miplib.zib.de/instance\underline{~}details\underline{~}scpn2.html}, except for one to be published in the next updates.}

\section{Theoretical Aspects of Our Algorithm}\label{theory}
In this section, we provide theoretical Aspects of our algorithm. While it is hard to have worst-case performance guarantee for general ILP heuristic algorithms, we can still reveal some underlying properties of our algorithm. 

{We are going to show the properties of our new operators and our algorithm on finite feasible ILP. Since our algorithm does not prove infeasibility or unboundedness, we do not analyze these two cases.} 
{In Algorithms 3-5 we can summarize that there are four ways to generate a solution in our algorithm: 

(1) apply a tight move operator to an infeasible solution 

(2) apply a lift move operator to a feasible solution 

(3) apply a unit incremental move when no positive $lm$ operation is found in \improve{} mode 

(4) move a variable's value to one side of its global bound

Following the observations in Section~\ref{convex}, we define the concept of boundary solutions, and show that all feasible solutions visited by our algorithm are boundary solutions. 
The proofs of all propositions are presented in the appendices.

\subsection{Boundary Solutions}

We assume an instance of ILP of the form: 
\begin{equation} \label{eq2}
\begin{split}
Minimize\ \ \    &obj(\bl{x})=\bl{c}^\top \bl{x} \\ 
subject\ to\ \ \ \    &\bl{A}\bl{x} \leq \bl{b} \\
                       &\bl{x} \in \mathbb{Z}^n  \\
\end{split}
\end{equation} 
For the convenience of analysis, here we do not distinguish variables' bounds from other constraints, {thus all variables' bounds are included in $\bl{A}\bl{x} \leq \bl{b}$, this is slightly different in Formula~(\ref{eq1}), and we take Formula~(\ref{eq2}) as the description of an ILP instance within Section~\ref{theory}.} 
We assume at least one of the coefficients in the objective function is non-zero: $\exists j \in \{1, ..., n\}$, $c_j\neq 0$, otherwise {this ILP instance has a constant objective function $0$ and turns to be a {satisfiability} problem. Since our algorithm does not prove feasibility, we exclude this case in analyzing our algorithm.} We denote the index set $I=\{1, ..., m\}$, $J=\{1, ..., n\}$, $J_{\neq 0}=\{ j | c_j \neq 0, j=1, ..., n\}$, $J_{\neq 0} \neq \emptyset$. Moreover, we assume that there is no variable that is free (does not appear in any constraint, has both infinite upper/lower bounds, and does not appear in $\bl{c^\top} \bl{x}$). 

We assume the ILP has a finite optimal solution $\bl{x}^*$ and a finite optimal objective value $opt$. We assume $\bl{x}^*$ exists (while it does not have to be unique), as our algorithm aims to find good solutions and does not try to prove an ILP is infeasible. {We assume there is a finite optimal objective value; otherwise there is no meaningful optimal objective value.} 

Now we present the concept of boundary solutions that will be used to analyze our operators, and show that all optimal solutions are boundary solutions. 

For an ILP instance in the form of Formula (\ref{eq2}), let polyhedron $P=\{ \bl{x} \in \mathbb{R}^n | \bl{A} \bl{x} \leq \bl{b}\}$, then the set of feasible solutions of the ILP could be described as $P \cap \mathbb{Z}^n$, and all feasible solutions belong to the \textbf{integer hull} $P_I=conv(P \cap \mathbb{Z}^n)$ where $conv(S) \subseteq \mathbb{R}^n$ for $S \subseteq \mathbb{Z}^n$ denotes the convex hull of a set of points $S$. Let $U = \cup \{\bl{e}_j, -\bl{e}_j\}$, $j\in J $. We call $\bl{x}+\bl{d}, \bl{d} \in U$ the \textbf{neighbors} of $\bl{x}$. Given $P=\{ \bl{x} \in \mathbb{R}^n | \bl{A} \bl{x} \leq \bl{b}\}$, we define the set of boundary points of $P$: 

\begin{myDef}
  $\bl{x} \in \mathbb{Z}^n$ is a \textbf{boundary point} of $P$ 
  if $\bl{x} \in P$ and  $\exists \bl{d} \in U$, $\bl{x}+\bl{d} \notin P$. The set of boundary points of $P$ is denoted by $\delta(P)$. 
\end{myDef}

This definition says a boundary point has at least one neighbor that is out of the feasible region, which is similar to the definition of boundary in topology.

\begin{figure}[ht]
    \centering
    \includegraphics[width=0.3\textwidth]{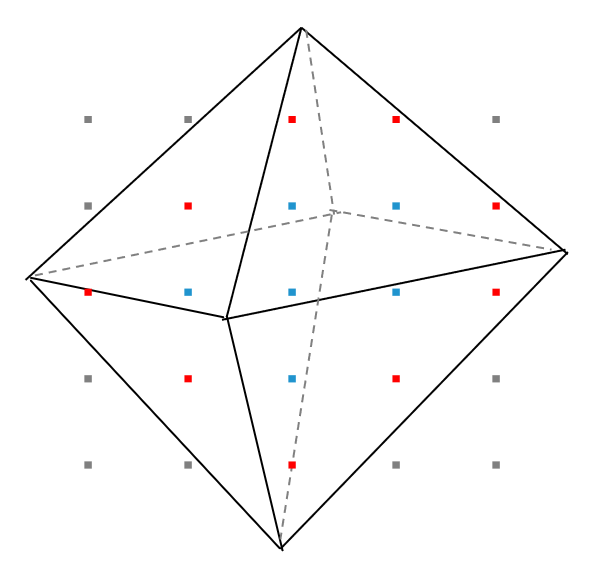}
    \caption{\centering{An illustration of a polyhedron of LP relaxation of an ILP and different solutions: grey/red/blue points represent infeasible/boundary/other solutions }}
    \label{fig: boundary}
\end{figure}

We call a solution $\bl{x}$ of Formula (\ref{eq2}) a \textbf{boundary solution} if and only if $\bl{x}$ is a boundary point of the polyhedron of its LP relaxation $P$. 
We can see in Figure \ref{fig: boundary} an illustration of different types (boundary/infeasible/others) of solutions. 
Then we exhibit the first property of ILP we use: every optimal solution is a boundary solution.

\begin{proposition}\label{prp1}
    For an ILP instance as the Formula (\ref{eq2}), its any optimal solution is a boundary solution, and also
    $\min_{\bl{x}\in P \cap \mathbb{Z}^n } obj(\bl{x}) =  \min_{\bl{x}\in \delta(P)\cap \mathbb{Z}^n } obj(\bl{x}) $. 
\end{proposition}

This proposition states that $\delta(P)$ is ``complete'' for obtaining an optimal solution: 

\begin{remark}
    The search space of an ILP instance in the form of Formula (\ref{eq2}) can be reduced to $\delta(P)$ and all optimal solutions are kept. Then $\delta(P)$ can be considered a \textbf{complete} space for ILP, in the sense that it would not miss the optimal solution. 
\end{remark}
We also show that the concept of $\delta(P)$ is significant: it is not a definition that is so general that it naturally contains all optimal solutions. We consider a polyhedral $P=\{ \bl{x} \in \mathbb{R}^n | \bl{A} \bl{x} \leq \bl{b}\}$, there could be different ILP instances specified by $P$ with different objective functions. We show a simple fact that for given $P$ and $\delta(P)$, any smaller subset of $\delta(P)$ may miss an optimal solution for some ILP instances consisting of $P$ with some objective function. 
\begin{fact}\label{mindelta}
    For polyhedral $P=\{ \bl{x} \in \mathbb{R}^n | \bl{A} \bl{x} \leq \bl{b}\}$, for any boundary point $\bl{x} \in \delta(P)$, there exist {a linear objective function $obj(\bl{x})$ such that $\bl{x}$ is the optimal solution of the ILP instance consisting of $P$ and this objective function.}   
\end{fact}
Since $\bl{x}\in P$ and $\bl{x}+\bl{e}_j \notin P$ for some $j$, it is easy to obtain such an instance by just setting the objective function $obj(\bl{x})=\bl{e}_j^\top \bl{x}$. It's trivial to see that $\bl{x}$ is optimal. 
{\begin{fact}\label{pbo}
    For a $0-1$ ILP instance (all variables are binary), all feasible solutions are boundary solutions.  
\end{fact}
}
\subsection{{Properties of our Operators}}
Now we show how our operators make use of this property. We first analyze \textit{tight move} operator and \textit{lift move} operator individually, and then show that all feasible solutions obtained by our algorithm are boundary solutions. 

\noindent \textbf{Tight Move Operator.} We show a property of the \textit{tight move} operator: if a solution obtained by the \textit{tight move} operator is feasible, then it is a boundary solution. 

We introduce some notions to facilitate the analysis of operators. Given $\bl{x} \in \mathbb{Z}^n$, $\forall j \in J, i \in I, A_{ij}\neq 0$, let $\phi_{ji}: \mathbb{Z}^n \rightarrow  \mathbb{Z}^n$ s.t. $\phi_{ji}(\bl{x})=\bl{x}'$ if $\bl{x}'$ is obtained from $\bl{x}$ by $tm(x_j, con_i, \bl{x}) $. 

\begin{proposition}\label{prp2}
     Any feasible solution obtained by \tm{} operator is a boundary solution: for $\bl{x} \in \mathbb{Z}^n$, $\forall j \in J, i \in I$, if $ \bl{x}'=\phi_{ji}(\bl{x})$ and $\bl{x}' \in P$, then $\bl{x}' \in \delta(P)$. 
\end{proposition}

\noindent \textbf{Lift Move Operator.} For the \textit{lift move} operator, we show another property, that it maps a feasible solution to a boundary solution: 

Given $\bl{x} \in \mathbb{Z}^n$, $\forall j \in J$, let $\chi_{j}: \mathbb{Z}^n \rightarrow  \mathbb{Z}^n$ s.t. $\chi_{j}(\bl{x})=\bl{x}'$ if $\bl{x}'$ is obtained from $\bl{x}$ by $lm(x_j, \bl{x})$.  
\begin{proposition}\label{prp3}
     The \textit{lift move} operator maps a feasible solution to a boundary solution: for $\bl{x} \in P\cap \mathbb{Z}^n$, $\forall j \in J$, $\chi_{j}(\bl{x}) \in \delta(P)$. 
\end{proposition}
\noindent \textbf{Our Algorithm.} As seen from Algorithms \ref{ls-3} - \ref{ls-4}, our algorithm adopts a strategy to apply the\textit{ lift move} operator for feasible solutions and \textit{tight move} for infeasible solutions. This coincides with the properties we showed in the preceding sections; additionally, we show here that the perturbations are also consistent with this strategy, and we can have the following argument for our algorithm:  
\begin{proposition}\label{algospace}
    Any feasible solution obtained in our algorithm is a boundary solution. 
\end{proposition}
This property shows our algorithm avoids visiting points in $P \setminus \delta(P)$, which is very helpful for problems with large variable domains. 

\subsection{Different Behaviors in Three Modes}

The three modes of our algorithm have different functionalities. The \search{} mode aims to find a feasible solution, i.e., locate $P$ in the full domain of variables; the \improve{} mode tries to generate a better solution while keeping the feasibility, i.e., walk  inside $P$ in the direction that improves the objective function; 
and the \restore{} mode wants to reach another feasible solution in $P$ again when an infeasible solution is reached from a perturbation in \improve{} mode. 
Thus, different settings of operators provide different behaviors in the three modes, adapted to different purposes. 

\subsubsection{Scoring Functions}

For \search{} and \restore{} modes, they both aim at finding feasible solutions, and both use \tm{} operator. The \tm{} scoring function is then focused on minimizing the degree of violation of constraints, which could be regarded as a ``distance'' to boundary solutions of $P$. For \improve{} mode, the purpose is to improve the objective function. Since solutions are all feasible in \improve{} mode, the only measure is the objective value, which is set as the \lm{} scoring function. 

\subsubsection{Step Size and Perturbations}
As our algorithm does not have a fixed step size to generate new solutions, in each mode, the actual step size and degree of perturbation depend on different operators. 

For \search{} and \restore{} modes which aim at the feasible solution, a larger discrepancy should be realized by the search strategy to get a rough location of $P$ quickly. This is handled in Algorithm \ref{ls-3} line 5 and Algorithm \ref{ls-4} line 7, allowing them to pick non-positive \tm{} operations. This provides a larger degree of perturbation and results in a larger discrepancy in the full domain. 

For \improve{} mode, since we showed the \lm{} operator always reaches boundary solutions, the step size is automatically adapted to the actual status. Once the \improve{} mode encounters a local optimum, i.e., no \lm{} operation to choose to stay in $P$, then it goes out of the region of $P$ by taking a unit incremental move in a variable with only step size $1$, to jump out of the local optimum. 
Since all feasible solutions are in $P$, the smallest step size is set here to not go out at a large distance from $P$.
Therefore, the smallest step size helps to enter $P$ again. 


\section{Conclusions and Future Work
\label{Conclusions}
}
This work proposed new characterizations of ILP with the concept of boundary solutions. 
Motivated by the new characterizations, we proposed a new local search algorithm Local-ILP, which is an efficient \PL{incomplete solver} for general integer linear programming validated on widely different types of problems. 
The main features of our algorithm include a new framework adopting three different modes, and two new operators designed for general ILP, the \textit{tight move} and \textit{lift move} operators, with tailored scoring functions for each. Experiments show that, in solving large-scale hard ILP problems within a reasonably short time, our algorithm is competitive and complementary to the state-of-the-art commercial solver Gurobi, and significantly outperforms the state-of-the-art non-commercial solver SCIP. More encouragingly, our algorithm established new records for 6 MIPLIB open instances.
We also presented the theoretical analysis of our algorithm, which shows our algorithm could avoid visiting unnecessary regions.

\PL{As for future work, we would like to design more sophisticated operators and scoring functions related to the constraints and the objective function, which we believe can further improve the performance of local search algorithms for ILP.
Given the promising results of Local-ILP, we would like to apply our local search framework to other combinatorial search problems, such as mixed integer programming, which is a more generalized model.}

\section*{Acknowledgement}
This work is supported by the Strategic Priority Research Program of the Chinese Academy of Sciences, Grant No. XDA0320000 and XDA0320300, and NSFC Grant 62122078.
\appendix

\section{Proof of Proposition}
\subsection{Proof of Proposition 1}
Let's consider an optimal solution $\bl{x}^*$ of Formula (2), as we assumed, $\bl{x}^*$ exists and is finite. For a $j\in J_{\neq 0}$, let $\bl{d}=\bl{e}_j$ if $c_j<0$, and $\bl{d}=-\bl{e}_i$ if $c_j>0$. 
Assume $\bl{x}^*+\bl{d} \in P$, then since $\bl{x}^*+\bl{d} \in \mathbb{Z}^n $, $\bl{x}^*+\bl{d}$ is a feasible solution of Formula (2) and $$obj(\bl{x}^*+\bl{d})=\bl{c}^T\bl{x}^*+\bl{c}^T\bl{d}=obj(\bl{x}^*)-|c_j|< obj(\bl{x}^*)$$ then there is another feasible solution with strictly smaller objective, contradiction with the assumption that $\bl{x}^*$ is optimal, so $\bl{x}^*+\bl{d} \notin P$ and $\bl{x}^*$ is a boundary solution. \qed 

\subsection{Proof of Proposition 2}
Let $\Delta=b_i- (\bl{A_i}) \cdot \bl{x}$,  
W.L.O.G. we assume $A_{ij}>0$, and $\Delta<0$ the other cases of $A_{ij}$ and $\Delta$ could be showed similarly. 

From Definition 1 we know $\bl{x}'=\bl{x}-min(\left| \left \lfloor {\Delta}/{A_{ij}}  \right \rfloor \right|, \left| x^l_j-x_j \right|) \cdot \bl{e}_j$. 
Since $A_{ij}>0$, we consider $\bl{x}'+\bl{e}_j$:

(1)If $min(\left| \left \lfloor {\Delta}/{A_{ij}}  \right \rfloor \right|, \left| x^l_j-x_j \right|)=\left| \left \lfloor {\Delta}/{A_{ij}}  \right \rfloor \right|$, 
since $\Delta<0, A_{ij}>0$, $\left| \left \lfloor {\Delta}/{A_{ij}}  \right \rfloor \right| = - \left \lfloor {\Delta}/{A_{ij}}  \right \rfloor $
then 
$$\bl{A_i} \cdot (\bl{x}'+\bl{e}_j) = \bl{A_i} \cdot \bl{x}+ \bl{A_i}(\left \lfloor {\Delta}/{A_{ij}}  \right \rfloor +1 )\bl{e}_j >  \bl{A_i} \cdot \bl{x}+ A_{ij}( {\Delta}/{A_{ij}} )=b_i$$
that is $\bl{A_i} \cdot (\bl{x}'+\bl{e}_j)>b_i$, which means $\bl{x}'+\bl{e}_j \notin P$, since assumed $\bl{x}' \in P$, then $\bl{x}'\in \delta(P)$. 

(2) If $min(\left| \left \lfloor {\Delta}/{A_{ij}}  \right \rfloor \right|, \left| x^l_j-x_j \right|)=\left| x^l_j-x_j \right|$, $\bl{x}'=\bl{x}- (\left| x^l_j-x_j \right|) \cdot \bl{e}_j$, we consider $\bl{x}'-\bl{e}_j$. Since $x_j$ always satisfy its global bound(by definition, all our algorithm will not break the global bound), $x^l_j-x_j \leq 0 $, then 
$$(\bl{x}'-\bl{e}_j)_j = x_j - \left| x^l_j-x_j \right| - 1 = x_j + (x^l_j-x_j) - 1 < x^l_j $$
$(\bl{x}'-\bl{e}_j)$ violated the global bound of the variable $x_j$, so $\bl{x}'-\bl{e}_j \notin P$, since assumed $\bl{x}' \in P$ then $\bl{x}'\in \delta(P)$. \qed

\subsection{Proof of Proposition 3}
Let $ \bl{x}'=\chi_{j}(\bl{x})$, $\Delta=b_i- (\bl{A_i}) \cdot \bl{x}$. W.L.O.G assume $c_j<0$ and $A_{ij}>0$,  from Definition 5, $lm(x_j, \bl{x})$ changes $\bl{x}'_j$ to the upper bound of $lfd(x_j,\bl{x})$. Since $lfd(x_j,\bl{x})=(\cap_i ldc(x_j,con_i,\bl{x}) ) \cap [x_j^l, x_j^u]$, the upper bound of $lfd(x_j,\bl{x})$ is either the upper bound of $[x_j^l, x_j^u]$ or of $ldc(x_j,con_i,\bl{x})$ for some $i$.

In the former case, $\bl{x}'+\bl{e}_j$ exceeds the upper bound of $[x_j^l, x_j^u]$. 

In the latter case,  from definition $ldc(x_j,con_i,\bl{x}) = \left( -\infty , x_j +  \left \lfloor {\Delta}/{A_{ij}} \right \rfloor \right] $ thus $\bl{x}' =  \bl{x} +  \left \lfloor {\Delta}/{A_{ij}} \right \rfloor \bl{e}_j$, then 
$$\bl{A_i}(\bl{x}'+\bl{e}_j)=\bl{A_i}(\bl{x} + \left \lfloor {\Delta}/{A_{ij}} \right \rfloor \bl{e}_j+\bl{e}_j) > \bl{A_i}(\bl{x} +  {\Delta}/{A_{ij}}  \bl{e}_j) = \bl{A_i} \bl{x} + \Delta = b_i $$
thus $\bl{A_i}(\bl{x}'+\bl{e}_j) > b_i$ which means $\bl{x}'+\bl{e}_j \notin P$. 

So in both case, we have $\bl{x}'+\bl{e}_j \notin P$. Moreover, it's easy to check $\bl{x}' \in P$ by definition that $ldc()$ is computed satisfying all constraints and global bound for the variable's bounds, and thus $\bl{x}' \in \delta(P) $. \qed

\subsection{Proof of Proposition 4}
In Algorithm 3 - 5 we can see that our algorithm has 4 ways to generate a new assignment:

(1) apply a tight move operator to an infeasible solution

(2) apply a lift move operator to a feasible solution

(3) (perturbation) apply a unit incremental move when no positive $lm$ operation is found in Improve mode

(4) (perturbation) move a variable's value to one side of its global bound

From Proposition 2 and Proposition 3, we know that the feasible solutions generated by case (1) and (2) must be boundary solutions. 

The case (3) will generate an infeasible solution because when there is no positive $lm$ operation, there is no operation to modify one variable's value to get a better feasible solution. In this case, a unit incremental move must generate an infeasible solution; otherwise it contradicts the condition that no positive $lm$ operation is found. 

In case (4), it is trivial that the generated solution is a boundary solution if it is feasible as one variable is set to one side of its global bound. 

In total, all feasible solutions generated by our algorithm are boundary solutions. 
\qed

 \bibliographystyle{elsarticle-num-names} 
 \bibliography{ref}





\end{document}